\documentclass[
a4paper,
11pt,
twoside,
]{article}
\usepackage[utf8]{inputenc}
\usepackage[T1]{fontenc}
\usepackage[osf,sc,noBBpl]{mathpazo} 
\usepackage[scaled]{beramono}
\usepackage{microtype}
\usepackage[english]{babel}
\usepackage{geometry}
\usepackage{amsmath,amsfonts,amssymb,amsthm}

\usepackage{etoolbox}

\usepackage{graphicx}
\usepackage[dvipsnames]{xcolor}

\theoremstyle{plain}

\usepackage{mathtools}
\usepackage{amssymb}
\usepackage{siunitx}
\usepackage{booktabs}

\newcommand{\D}{\mathop{}\!\mathrm{d}}

\newcommand{\I}{\mathrm{i}}
\newcommand{\E}{\mathrm{e}}

\renewcommand{\restriction}{\raise-.5ex\hbox{\ensuremath{\upharpoonright}}}

\definecolor{webgreen}{rgb}{0,.5,0}
\definecolor{webbrown}{rgb}{.6,0,0}
\definecolor{myblue}{rgb}{0,0.25,0.5}

\usepackage[
  colorlinks=true,
  breaklinks=true,
  bookmarksnumbered,
  pdfhighlight=/O,
  urlcolor=webbrown,
  linkcolor=myblue,
  citecolor=webgreen,
  hyperfootnotes=true,
]{hyperref}
\newcommand{\email}[1]{\href{mailto:#1}{\texttt{#1}}}
\makeatletter
\newcommand{\pagerefstar}{\@pagerefstar}
\makeatother

\usepackage{fnpct}

\usepackage{xcolor}
\definecolor{lightergray}{gray}{0.97}
\definecolor{darkergray}{gray}{0.4}
\usepackage{listings}
\usepackage{csquotes}

\usepackage[
  mincrossrefs=99,
  style=numeric-comp,
  bibstyle=siam,
  sorting=nyt,
  backref=true,
  backend=biber,
  url=true,
  doi=true,
  maxbibnames=99,
  maxnames=99,
  giveninits=true,
  useprefix=true,
  date=year,
  autolang=other
]{biblatex}
\DeclareCiteCommand{\citeauthornsc}
  {%
   \boolfalse{citetracker}%
   \boolfalse{pagetracker}%
   \usebibmacro{prenote}}
  {\ifciteindex
     {\indexnames{labelname}}
     {}%
   \printnames{labelname}}
  {\multicitedelim}
  {\usebibmacro{postnote}}
\DeclareSourcemap{
  \maps[datatype=bibtex]{
    \map{
      \step[fieldsource=doi,final]
      \step[fieldset=url,null]
    }
  }
}

\usepackage{caption}
\usepackage{fancyhdr}
\usepackage{lastpage}

\title{A practical guide to piecewise pseudospectral collocation for Floquet multipliers of delay equations in MATLAB}

\author{Dimitri Breda, Davide Liessi and Rossana Vermiglio}

\newcommand{\mydate}{22 March 2022}

\fancypagestyle{mypagestyle}{%
  \fancyhf{}%
  \fancyhead[LO,RE]{\scriptsize Practical guide to piecewise pseudospectral collocation for Floquet multipliers}%
  \fancyhead[RO,LE]{\scriptsize Breda, Liessi, Vermiglio}%
  \fancyfoot[LE,RO]{\scriptsize \thepage\,/\,\pagerefstar{LastPage}}%
  \fancyfoot[LO,RE]{\scriptsize \mydate{}}%
}
\begin{document}

\thispagestyle{empty}

\begin{center}
\LARGE
A practical guide to piecewise pseudospectral collocation for Floquet multipliers of delay equations in MATLAB

\bigskip
\large
Dimitri Breda\footnote{\email{dimitri.breda@uniud.it}},
Davide Liessi\footnote{\email{davide.liessi@uniud.it}},
Rossana Vermiglio\footnote{\email{rossana.vermiglio@uniud.it}}

\medskip
\small
CDLab -- Computational Dynamics Laboratory \\
Department of Mathematics, Computer Science and Physics, University of Udine \\
Via delle Scienze 206, 33100 Udine, Italy

\medskip
\large
\mydate
\end{center}

\begin{abstract}
In recent years we provided numerical methods based on pseudospectral collocation for computing the Floquet multipliers of different types of delay equations, with the goal of studying the stability of their periodic solutions.
The latest work of the series concerns the extension of these methods to a piecewise approach, in order to take the properties of numerically computed solutions into account.
In this chapter we describe the MATLAB implementation of this method and provide practical usage examples.
\end{abstract}

\section{Introduction}

Delay equations appear often in mathematical modeling: indeed, the presence of delay terms allows to base the current evolution on the past history, increasing the realism of the model.
Examples of fields where delays arise naturally are control theory in engineering, e.g.
\cite{%
Fridman2014,%
GuKharitonovChen2003,%
MichielsNiculescu2014,%
Stepan1989%
},
and population dynamics or epidemics in mathematical biology, e.g.
\cite{%
ArinoVanDenDriessche2006,%
BredaDiekmannDeGraafPuglieseVermiglio2012,%
Kuang1993,%
MacDonald1978,%
MetzDiekmann1986,%
Smith2011%
}
(see
\cite{%
Erneux2009,%
KolmanovskiiMyshkis1999%
}
for further applications).

These models often present periodic behaviors, whose emergence is particularly facilitated by the presence of delays, and there is a strong interest in applications to studying their asymptotic stability.
Examples include problems of network consensus, mechanical vibrations, endemic states and seasonal fluctuations.

In the case of delay differential equations (DDEs), the main numerical tool is DDE-BIFTOOL%
\footnote{\url{http://ddebiftool.sourceforge.net/}}
\cite{%
EngelborghsLuzyaninaRoose2002,%
SieberEngelborghsLuzyaninaSamaeyRoose2014%
},
which bases its methods for periodic solutions on the piecewise orthogonal collocation of \cite{EngelborghsLuzyaninaIntHoutRoose2001}.
Other approaches to the stability of periodic solutions include the semi-discretization method \cite{InspergerStepan2011} and Chebyshev-based collocations \cite{Bueler2007,ButcherBobrenkov2011,ButcherMaBuelerAverinaSzabo2004}, and more references can be found in \cite{Jarlebring2008}, while a comparison can be found in \cite{LehotzkyInsperger2016}.
Particularly relevant for the present work is \cite{BredaMasetVermiglio2012}, based on the pseudospectral collocation of evolution operators, which is possibly one of the most general in terms of combinations of discrete and distributed delays.
Moreover, the recent work \cite{BorgioliHajduInspergerStepanMichiels2020} provides a generalization of the collocation approach of DDE-BIFTOOL, allowing for discontinuous coefficients.

Extensions of \cite{BredaMasetVermiglio2012} have been recently proposed for renewal equations (REs) \cite{BredaLiessi2018} and systems of REs and DDEs \cite{BredaLiessi2020}.
In fact, as far as we know, these are the first available methods for REs.
In \cite{BredaLiessiVermiglio} we further extended these methods, which are based on polynomials on a single piece, to the case of piecewise polynomials, in order to take into account the partition of the period interval used for the numerically approximated periodic solutions: therein we show that this is essential if the said partition is far from uniform, as observed also in \cite{BorgioliHajduInspergerStepanMichiels2020}.

In this work we present a MATLAB implementation of \cite{BredaLiessiVermiglio}, available at \url{http://cdlab.uniud.it/software}, and provide some practical usage examples, showcasing also its generality with respect to the type of equations (DDEs, REs, systems of REs and DDEs) and to the application to studying the dynamics (stability of equilibria and periodic solutions, detection of bifurcations and stability boundaries).
After presenting the basics theoretical aspects of the problem in Sect.~\ref{ch:BLV:sec:theory}, we describe the main aspects of the reformulation of the evolution operator and its discretization in Sect.~\ref{ch:BLV:sec:discr}.
Then in Sect.~\ref{ch:BLV:sec:impl} we present the MATLAB implementation and describe some of our choices and of the user-available options.
Finally, Sect.~\ref{ch:BLV:sec:ex} presents the usage examples with usable code excerpts, followed by some comments on advanced use cases in Sect.~\ref{ch:BLV:sec:adv}.

\section{Theoretical foundations}
\label{ch:BLV:sec:theory}

Let $d_{X}$ and $d_{Y}$ be nonnegative integers, not both null, $\tau$ a positive real.
We consider systems of REs and DDEs
\begin{equation}
\label{ch:BLV:eq:nl}
\left\{
\begin{aligned}
&x(t) = F(x_{t}, y_{t}), \\
&y'(t) = G(x_{t}, y_{t}),
\end{aligned}
\right.
\end{equation}
for
$F \colon X \times Y \to \mathbb{R}^{d_{X}}$,
$G \colon X \times Y \to \mathbb{R}^{d_{Y}}$,
$X \coloneqq L^{1}([-\tau, 0], \mathbb{R}^{d_{X}})$ with the $L^{1}$ norm,
$Y \coloneqq C([-\tau, 0], \mathbb{R}^{d_{Y}})$ with the uniform norm,
$x_{t}(\theta) \coloneqq x(t+\theta)$ for $\theta \in [-\tau, 0]$
and $y_{t}$ defined likewise.
In the following we call \eqref{ch:BLV:eq:nl} a \emph{coupled equation}.
Observe that \eqref{ch:BLV:eq:nl} reduces to a DDE if $d_{X} = 0$ and to an RE if $d_{Y} = 0$.

Assume that \eqref{ch:BLV:eq:nl} has an $\omega$-periodic solution $(\bar{x}, \bar{y})$ with corresponding linearization
\begin{equation}
\label{ch:BLV:eq:lz}
\left\{
\begin{aligned}
&x(t) = DF(\bar{x}_{t}, \bar{y}_{t})(x_{t}, y_{t}), \\
&y'(t) = DG(\bar{x}_{t}, \bar{y}_{t})(x_{t}, y_{t}),
\end{aligned}
\right.
\end{equation}
where $D$ indicates the Fréchet derivative%
\footnote{%
For a discussion of the difficulties of linearizing REs, see \cite[Sect.~3.5]{DiekmannGettoGyllenberg2008} about equilibria; extending those results to periodic solutions is an open problem.
}%
.
Observe that \eqref{ch:BLV:eq:lz} is an $\omega$-periodic nonautonomous equation.
The information on the local stability of the periodic solution can be derived from the spectrum of the associated monodromy operators, whose nonzero elements are called Floquet multipliers (\emph{multipliers} in the following).
The monodromy operators are given by $U(s+\omega, s)$, where $U(t,s) \colon X \times Y \to X \times Y$ for $t \geq s$ is the evolution operator associated to \eqref{ch:BLV:eq:lz}, defined as
\begin{equation*}
U(t, s)(\phi, \psi) = (x(\cdot; s, (\phi, \psi))_{t}, y(\cdot; s, (\phi, \psi))_{t}),
\end{equation*}
for $(x(\cdot; s, (\phi, \psi)), y(\cdot; s, (\phi, \psi)))$ the solution of the initial value problem for \eqref{ch:BLV:eq:lz} with initial condition $(x_{s}, y_{s}) = (\phi, \psi)$%
\footnote{%
For the well-posedness of the initial value problem for coupled equations, see \cite[Theorem 1]{BredaLiessi2020}.
}%
.

For DDEs, assuming that $G$ is a $C^{1}$ function, the multipliers are independent of $s$ \cite[Theorem XIII.3.3]{DiekmannVanGilsVerduynLunelWalther1995} and $1$ is always a multiplier (called \emph{trivial}) due to the linearization around a periodic solution \cite[Theorem XIV.2.6]{DiekmannVanGilsVerduynLunelWalther1995}.
Thanks to the principle of linearized stability, provided that the trivial multiplier is simple, the periodic solution is locally asymptotically stable if all nontrivial multipliers are inside the unit circle, while it is unstable if at least one multiplier is outside \cite[Sect.~XIV.3 and Exercise XIII.2.3]{DiekmannVanGilsVerduynLunelWalther1995}.
The same results hold for REs assuming that $F$ is a $C^{1}$ and globally Lipschitz continuous function \cite{BredaLiessi2018,BredaLiessi2021}.
For coupled equations a Floquet theory is currently missing, but it is reasonable to expect that similar results are valid, as we assume henceforth.

\section{Discretization}
\label{ch:BLV:sec:discr}

The following reformulation and discretization can be applied to any evolution operator; however, they are typically applied to monodromy operators, hence our use of $\omega$ in the sequel, even if it need not be the period.
Let $T \coloneqq U(s+\omega, s)$ be an evolution operator for $s \in \mathbb{R}$ and $\omega \geq 0$.
In order to obtain a finite-dimensional approximation of $T$, we first reformulate the operator as follows.
Let
$X^{+} \coloneqq L^{1}([0, \omega], \mathbb{R}^{d_{Y}})$
and
$X^{\pm} \coloneqq L^{1}([-\tau, \omega], \mathbb{R}^{d_{Y}})$
with the corresponding $L^{1}$ norms, and
$Y^{+} \coloneqq C([0, \omega], \mathbb{R}^{d_{Y}})$
and
$Y^{\pm} \coloneqq C([-\tau, \omega], \mathbb{R}^{d_{Y}})$
with the corresponding uniform norms.
Let $V\colon (X \times Y) \times (X^{+} \times Y^{+}) \to X^{\pm} \times Y^{\pm}$ be the operator which, given pairs of functions $(\phi, \psi)$ on $[-\tau, 0]$ (think of it as the initial value) and $(w, z)$ on $[0, \omega]$ (think of it as the result of applying the right-hand side of \eqref{ch:BLV:eq:lz} to the solution), constructs the solution of \eqref{ch:BLV:eq:lz} on $[-\tau, \omega]$ as
\begin{equation}
\label{ch:BLV:eq:V}
V((\phi, \psi), (w, z))(t) \coloneqq
\begin{cases}
\left(w(t), \displaystyle \psi(0) + \int_{0}^{t} z(\sigma) \D\sigma\right), & t \in (0, \omega], \\
(\phi(t), \psi(t)), & t \in [-\tau, 0].
\end{cases}
\end{equation}
Let also $\mathcal{F}_{s} \colon Y^{\pm} \times Y^{\pm} \to X^{+} \times Y^{+}$ be the operator defined as
\begin{equation*}
(\mathcal{F}_{s}(u, v)) (t) \coloneqq (DF(\bar{x}_{s+t}, \bar{y}_{s+t})(u_{t}, v_{t}), DG(\bar{x}_{s+t}, \bar{y}_{s+t})(u_{t}, v_{t})), \quad t\in[0,\omega],
\end{equation*}
which basically applies to its argument the action of the right-hand side of \eqref{ch:BLV:eq:lz} (with the time shifted by $s$ so that the initial time is $0$).
Finally, $T$ can be reformulated as
\begin{equation}
\label{ch:BLV:eq:T-as-V}
T (\phi, \psi) = V((\phi, \psi), (w^{\ast}, z^{\ast}))_{\omega},
\end{equation}
where $(w^{\ast}, z^{\ast}) \in X^{+} \times Y^{+}$ is the solution of the fixed point equation
\begin{equation}
\label{ch:BLV:eq:fixed-point}
(w, z) = \mathcal{F}_{s} V((\phi, \psi), (w, z)),
\end{equation}
which exists and is unique in the same conditions as the solutions of \eqref{ch:BLV:eq:lz}.
Observe that $w^{\ast}$ and $z^{\ast}$ are, respectively, the first component and the derivative of the second component of the solution of \eqref{ch:BLV:eq:lz} with initial functions $(x_{s}, y_{s})=(\phi, \psi)$.

\bigskip

The discrete counterparts of functions in $X$ and $Y$ are chosen to be real vectors of certain dimensions, depending both on the dimensions $d_{X}$ and $d_{Y}$ of \eqref{ch:BLV:eq:lz} and on some discretization parameters: let $d_{X}^{\bullet}$ and $d_{Y}^{\bullet}$ be those dimensions, respectively, so that $X$ is discretized as $X^{\bullet} \coloneqq \mathbb{R}^{d_{X}^{\bullet}}$ and $Y$ as $Y^{\bullet} \coloneqq \mathbb{R}^{d_{Y}^{\bullet}}$.
Let $R \colon X \times Y \to X^{\bullet} \times Y^{\bullet}$ be the restriction operator associating a pair of functions with their discretizations and $P \colon X^{\bullet} \times Y^{\bullet} \to X \times Y$ be the prolongation operator associating a pair of vectors to a pair of functions in a certain subset of $X \times Y$, typically a subspace of (piecewise) polynomials, such that $RP$ is the identity on $X^{\bullet} \times Y^{\bullet}$.
The operator $\mathcal{L} = PR$ is the approximation operator on $X \times Y$ defined by the chosen discretization scheme.
Similarly, discretizations $X^{+\bullet}$ and $Y^{+\bullet}$ of $X^{+}$ and $Y^{+}$ and corresponding operators $R^{+}$, $P^{+}$ and $\mathcal{L}^{+}$ are introduced.

Following \eqref{ch:BLV:eq:T-as-V} and \eqref{ch:BLV:eq:fixed-point}, the discretization of $T$ is the finite-dimensional operator $T^{\bullet} \colon X^{\bullet} \times Y^{\bullet} \to X^{\bullet} \times Y^{\bullet}$ defined as
\begin{equation}
\label{ch:BLV:eq:discrete-T-as-V}
T^{\bullet} (\Phi, \Psi) \coloneqq R V(P(\Phi, \Psi), P^{+}(W^{\ast}, Z^{\ast}))_{\omega},
\end{equation}
where $(W^{\ast}, Z^{\ast}) \in X^{+\bullet} \times Y^{+\bullet}$ is the solution of the fixed point equation%
\footnote{%
The well-posedness of \eqref{ch:BLV:eq:discrete-fixed-point} depends on the chosen discretization and needs to be proved.
}
\begin{equation}
\label{ch:BLV:eq:discrete-fixed-point}
(W, Z) = R^{+} \mathcal{F}_{s} V(P(\Phi, \Psi), P^{+}(W, Z))
\end{equation}
for the given $(\Phi, \Psi) \in X^{\bullet} \times Y^{\bullet}$.
The eigenvalues of $T^{\bullet}$ are then computed with standard methods and considered as approximations of the multipliers.

\bigskip

In this work we follow a pseudospectral collocation approach, which means that the restriction operator evaluates the functions at a fixed set of nodes, while the prolongation operator associates to these values the corresponding interpolating polynomial.

The version based on polynomials on a single piece was presented for DDEs in \cite{BredaMasetVermiglio2012,BredaMasetVermiglio2015}, for REs in \cite{BredaLiessi2018} and for coupled equations in \cite{BredaLiessi2020}, along with the proofs of the well-posedness of \eqref{ch:BLV:eq:discrete-fixed-point} and of the convergence of the approximated eigenvalues%
\footnote{%
In these works it is made clear that the method requires to evaluate only functions that are pointwise defined, even if the state space is $L^{1}$ for REs.
}%
.

In this work we present an implementation of the version based on continuous piecewise polynomials.
This choice is motivated by the fact that periodic solutions are usually computed as piecewise polynomials \cite{Ando2021,AndoBreda,AndoBreda2020b,Bader1985,EngelborghsLuzyaninaIntHoutRoose2001} and including the piece endpoints in the grid of discretization nodes is essential in some cases, as shown in \cite{BredaLiessiVermiglio}, where this discretization scheme is presented along with several examples and a discussion of the convergence properties (see Sect.~\ref{ch:BLV:sec:method} below for more details).

Under typical smoothness assumptions, all these methods exhibit spectral accuracy \cite{Trefethen2000}, with the order of convergence being infinite in the degree of the polynomials; moreover, the piecewise method exhibits a finite order of convergence in the number of pieces.

\section{Implementation}
\label{ch:BLV:sec:impl}

We provide a MATLAB implementation of the method described in \cite{BredaLiessiVermiglio}, which is available at \url{http://cdlab.uniud.it/software}.
In Sect.~\ref{ch:BLV:sec:equation} we describe our prototype linear coupled equation and how it is defined as an input to the method; then, in Sect.~\ref{ch:BLV:sec:method} we describe some details of the numerical method and how to control them.

\subsection{Defining the equation}
\label{ch:BLV:sec:equation}

As prototype problem we choose the linear coupled equation
\begin{equation}
\label{ch:BLV:eq:proto}
\left\{
\begin{aligned}
x(t) {}={}
&A_{XX}(t) x(t) + \sum_{k = 1}^{p} B^{(k)}_{XX}(t) x(t-\tau_{k}) + \sum_{k = 1}^{p} \int_{-\tau_{k}}^{-\tau_{k-1}} C^{(k)}_{XX}(t, \theta) x(t+\theta) \D \theta \\
&+ A_{XY}(t) y(t) + \sum_{k = 1}^{p} B^{(k)}_{XY}(t) y(t-\tau_{k}) + \sum_{k = 1}^{p} \int_{-\tau_{k}}^{-\tau_{k-1}} C^{(k)}_{XY}(t, \theta) y(t + \theta) \D\theta, \\
y'(t) {}={}
&A_{YX}(t) x(t) + \sum_{k = 1}^{p} B^{(k)}_{YX}(t) x(t-\tau_{k}) + \sum_{k = 1}^{p} \int_{-\tau_{k}}^{-\tau_{k-1}} C^{(k)}_{YX}(t, \theta) x(t+\theta) \D \theta \\
&+ A_{YY}(t) y(t) + \sum_{k = 1}^{p} B^{(k)}_{YY}(t) y(t-\tau_{k}) + \sum_{k = 1}^{p} \int_{-\tau_{k}}^{-\tau_{k-1}} C^{(k)}_{YY}(t, \theta) y(t + \theta) \D\theta,
\end{aligned}
\right.
\end{equation}
where $\tau_{0} \coloneqq 0 < \tau_{1} < \dots < \tau_{p} \coloneqq \tau$ are the delays, $A_{\square\triangle}, B^{(k)}_{\square\triangle} \colon \mathbb{R} \to \mathbb{R}^{d_{\square} \times d_{\triangle}}$ and $C^{(k)}_{\square\triangle} \colon \mathbb{R} \times [-\tau, 0] \to \mathbb{R}^{d_{\square} \times d_{\triangle}}$, with $\square$ and $\triangle$ standing for ``$X$'' or ``$Y$''.

Note that \eqref{ch:BLV:eq:proto} is actually of a more general type than \eqref{ch:BLV:eq:lz}%
\footnote{%
Indeed, thanks to the Riesz representation theorem for $L^{1}$ (see, e.g., \cite[p.~400]{RoydenFitzpatrick2010}), linear REs are integral equations with $L^{\infty}$ kernels and do not have current time or discrete delay terms (see also \cite[end of Sect.~2]{BredaLiessi2018}).
}%
.
However, from the point of view of the implementation, treating the terms with $A_{\square X}(t)$ and $B^{(k)}_{\square X}(t)$ is no more difficult than treating the others, while considering them in \eqref{ch:BLV:eq:proto} has the advantage of allowing the application of the method to \emph{neutral} REs, which are the object of ongoing research (see \cite{BredaLiessiVermiglio} for a first example).

In our MATLAB implementation the system of equations is specified as a structure array containing the system parameters as fields.
To help in checking the correctness of inputs and setting up some helper variables, we provide the \texttt{eigTMNpw\_system} function.
It takes four required positional parameters, namely the dimension of the system as the vector $\texttt{[}d_{X}\texttt{,}d_{Y}\texttt{]}$, the vector of delays in ascending order $\texttt{[}\tau_{1}\texttt{,\dots,}\tau_{p}\texttt{]}$, the length of time evolution $\omega$ and the initial time $s$.
It takes also several optional named parameters, corresponding to the coefficients (function handles \texttt{AXX}, \texttt{AXY}, \texttt{AYX}, \texttt{AYY}; cell arrays \texttt{BXX}, \texttt{BXY}, \texttt{BYX}, \texttt{BYY}, \texttt{CXX}, \texttt{CXY}, \texttt{CYX}, \texttt{CYY} of function handles, whose components are ordered according to the delays), the parameters (\texttt{par}, which can be a variable of any kind at the user's choice, depending on how it is used in the coefficients) and the partition of $[0, \omega]$ defining the piecewise approach (either the vector \texttt{t} of the partition endpoints or the number \texttt{L} of uniform pieces; more details in Sect.~\ref{ch:BLV:sec:method} below).
The \texttt{help} message for \texttt{eigTMNpw\_system} describes both the inputs and the output structure array more precisely.
We give examples in Sect.~\ref{ch:BLV:sec:ex}.

\subsection{Controlling the method}
\label{ch:BLV:sec:method}

The piecewise approach decribed in \cite{BredaLiessiVermiglio} and implemented here is based on a partition of the interval $[0, \omega]$, typically given by a numerical periodic solution of the equation computed as a piecewise polynomial.
Based on that partition, the interval $[-\tau, 0]$ is partitioned by subtracting $\omega$ (or multiples of $\omega$, if $\tau > \omega$) to the points in $[0, \omega]$.

\begin{figure}[htp]
\centering
\includegraphics{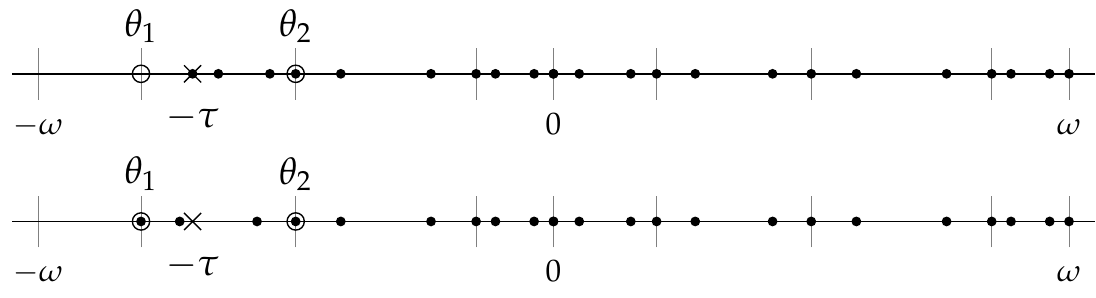}
\caption{%
Example partitions and collocation nodes with $\omega > \tau$.
Ticks mark the partition of $[0, \omega]$ and their translation by $-\omega$, crosses mark $-\tau$, circles mark $\theta_{1}$ and $\theta_{2}$, dots mark the collocation nodes.
The endpoint of the leftmost interval is chosen to be $-\tau$ (top) or $\theta_{1}$ (bottom).
}
\label{ch:BLV:fig:mesh}
\end{figure}

If $\tau$ does not coincide with one of the resulting points, a choice must be made on how to treat the leftmost incomplete piece of $[-\tau, 0]$.
In our implementation we provide two main approaches.
If $\theta_{1}$ and $\theta_{2}$ are the two endpoints of the provisional partition closest to $-\tau$ such that $\theta_{1} < -\tau < \theta_{2}$, the first approach is to consider the interval $[-\tau, \theta_{2}]$, while the second is to consider $[\theta_{1}, \theta_{2}]$, i.e., to artificially increase the maximum delay to $-\theta_{1}$%
\footnote{%
Observe that such an extension does not affect the dynamics.
}%
.
The latter may seem pointless, but it leads to a simpler expression of some of the matrices used to construct $T^{\bullet}$.
In our experiments we tested both approaches, obtaining (almost) identical results.
Figure~\ref{ch:BLV:fig:mesh} illustrates the two options.
Since the first approach may result in a very small interval, we provide also two variations of the first approach, differing for the treatment of the leftmost piece when its length is below a certain user-controllable threshold (more details in the \texttt{help} message for the \texttt{eigTMNpw\_method} function described below).
In the following, let $\tilde{\tau}$ be defined as $\tau$ or $\theta_{1}$ according to the chosen approach.

In each piece of the interval $[-\tilde{\tau}, \omega]$, we choose collocation nodes of the same family.
In order to approximate the operator $T$ and not only its spectrum, we need to impose the continuity conditions on the piecewise polynomials implicitly: the consequence is that the nodes need to include the interval endpoints.
The default choice in our implementation is Chebyshev extrema, but we provide also equidistant points and nodes based on Chebyshev and Gauss--Legendre zeros with the addition of the endpoints.

The numerical method is implemented in the \texttt{eigTMNpw} function, taking the structure array defining the equation described in Sect.~\ref{ch:BLV:sec:equation} above as required first positional parameter.
Other optional inputs are a positional and named parameter%
\footnote{%
Optional positional parameters can be specified also as named parameters in MATLAB but not in GNU Octave.
}
\texttt{method} controlling the details of the discretization and a named parameter \texttt{verbosity} to control the amount of information printed during the computation (see \texttt{eigTMNpw}'s \texttt{help} message).
The function outputs, in order, the vector of the eigenvalues of $T^{\bullet}$, $T^{\bullet}$ itself and four other matrices as described in Sect.~\ref{ch:BLV:sec:adv}.

In order to assist in correctly constructing the structure array \texttt{method}, we provide the \texttt{eigTMNpw\_method} function, whose input parameters are all optional.
Its first positional and named parameter \texttt{M} specifies the degree of the approximating piecewise polynomials, with a default value of $10$.
The choice of the family of collocation nodes is controlled by the named parameter \texttt{CollocationFamily}, while the strategy for the treatment of the leftmost incomplete piece of $[-\tilde{\tau}, 0]$ and, if needed, the threshold are chosen with the named parameters \texttt{Strategy} and \texttt{Threshold}, respectively.
The \texttt{Threshold} parameter (defaulting to \texttt{eps}) is used also to check that $\tau$ is not too small relative to $\omega$.

The parameters \texttt{QuadratureMethod} and \texttt{QuadratureParameters} control the choice of quadrature method and the relevant parameters: the default method is the Clenshaw--Curtis quadrature \cite{ClenshawCurtis1960,Trefethen2008}; other ready to use methods are MATLAB's \texttt{quad} and \texttt{integral}.

The last named parameter is \texttt{ZeroDirection}.
Indeed the formulation of the method requires to choose whether $0$ is included in the delay interval (i.e., intervals $[-\tau, 0]$ and $(0, \omega]$ are considered) or in the time evolution interval (i.e., intervals $[-\tau, 0)$ and $[0, \omega]$ are considered).
In all our experiments with DDEs, REs and coupled equations there was no difference between the two options and the first one would be in accordance with the formulation \eqref{ch:BLV:eq:V}.
However, when dealing with neutral REs the choice influences the resulting spectrum and the second option seems to provide more accurate results, motivating the fact that it is the default one.

The \texttt{help} message for \texttt{eigTMNpw\_method} describes both the inputs and the output structure array in more detail.

\subsection{Current limitations}

The coefficients of the equation are currently specified as matrix-valued functions, according to \eqref{ch:BLV:eq:proto}.
This prevents the quadrature method from using a vectorized approach when the equations are not scalar (i.e., either $d_{X}$ or $d_{Y}$ is greater than $1$).
This has a large impact on the efficiency of the code, since when integral terms are present, the quadrature needs to be performed several times and in multiple levels of integration.

Sometimes delays depend on the parameters of the system, so their order can change as the parameters vary.
While the code \texttt{eigTMN} for DDEs \cite{BredaMasetVermiglio2015} could automatically sort the delays and the corresponding coefficients, \texttt{eigTMNpw} is not yet capable of that.
As a result the user is forced to redefine the system each time the order of the delays changes.

Finally, there are several parts of the computation that may be further optimized, taking advantage of the specific form of the matrices used to construct $T^{\bullet}$ and reformulating some cycles to possibly spare some repeated quadratures.

\section{Usage examples}
\label{ch:BLV:sec:ex}

We now propose some examples of the use of \texttt{eigTMNpw}, showcasing its capabilities.
In Sect.~\ref{ch:BLV:sec:quadraticre} we compute the multipliers for an RE with a distributed delay term and an explicitly known periodic solution (hence not requiring the piecewise approach), and then find some bifurcations.
In Sect.~\ref{ch:BLV:sec:belairmackey} we use the piecewise approach for a DDE with two discrete delays and a periodic solution computed via DDE-BIFTOOL.
The example in \ref{ch:BLV:sec:logisticdaphniadelayed} is a coupled equation with terms of all the supported kinds (current time, discrete delays, distributed delays) and a periodic solution computed via MatCont and \cite{BredaDiekmannGyllenbergScarabelVermiglio2016}.
Finally, we compute in Sect.~\ref{ch:BLV:sec:hayes} the classic stability chart for the null equlibrium of the Hayes DDE \cite{Hayes1950}.
In all examples we used MATLAB R2019a.

\subsection{A renewal equation with an explicit periodic solution}
\label{ch:BLV:sec:quadraticre}

\begin{figure}[htp]
\centering
\includegraphics[trim=15.941pt 0 0 0]{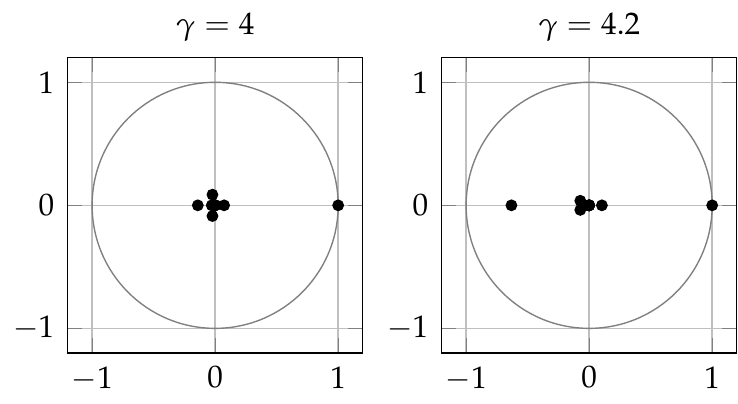}
\caption{%
Multipliers relevant to \eqref{ch:BLV:eq:quadraticre}--\eqref{ch:BLV:eq:quadraticresolution} computed as described in Sect.~\ref{ch:BLV:sec:quadraticre}.
}
\label{ch:BLV:fig:quadraticre}
\end{figure}

For our first example we consider the RE with quadratic nonlinearity
\begin{equation}
\label{ch:BLV:eq:quadraticre}
x(t) = \frac{\gamma}{2} \int_{-3}^{-1} x(t+\theta) (1-x(t+\theta)) \D\theta,
\end{equation}
which has a branch of $4$-periodic solutions with the explicit expression
\begin{equation}
\label{ch:BLV:eq:quadraticresolution}
\bar{x}(t) = \frac{1}{2} + \frac{\pi}{4 \gamma} + \sqrt{\frac{1}{2} - \frac{1}{\gamma} - \frac{\pi}{2 \gamma^{2}}\left(1+\frac{\pi}{4}\right)} \sin\left(\frac{\pi}{2} t\right),
\end{equation}
as proved in \cite{BredaDiekmannLiessiScarabel2016}.
To study the stability of $\bar{x}$, we consider the linear RE%
\begin{equation*}
x(t) = \frac{\gamma}{2} \int_{-3}^{-1} (1-2\bar{x}(t+\theta)) x(t+\theta) \D\theta.
\end{equation*}
We compute the multipliers relevant to \eqref{ch:BLV:eq:quadraticresolution} for $\gamma = 4$ with default options for the method with the following instructions.
\begin{lstlisting}
xbar = @(t, gamma) 1/2 + pi/(4*gamma) ...
  + sqrt(1/2 - 1/gamma - pi/(2*gamma^2)*(1+pi/4)) * sin(pi/2*t);
system = eigTMNpw_system(...
  [1, 0], ... % dimensions [dX, dY]
  [1, 3], ... % delays
  4, ...      % period (or length of time evolution)
  0, ...      % initial time
  'CXX', {[], @(t, theta, par) par/2 * (1 - 2*xbar(t+theta,par))}, ...
  'par', 4);
mult1 = eigTMNpw(system);
\end{lstlisting}
To compute the multipliers for another value of $\gamma$, the value of the parameter can be changed by setting \texttt{system.par} directly, without calling \texttt{eigTMNpw\_system} again.
\begin{lstlisting}
system.par = 4.2;
mult2 = eigTMNpw(system);
\end{lstlisting}
Figure~\ref{ch:BLV:fig:quadraticre} shows the multipliers computed as described above.
For $\gamma = 4$ the dominant ones are $1.000000179842839$, $-0.140831131942336$ and $-0.021890537332049\pm0.086918211021300\I$, while for $\gamma = 4.2$ they are $1.000000266174309$, $-0.631694832535750$ and $0.103689337250279$.

We can find some bifurcations by using MATLAB's \texttt{fzero} on a function testing for one multiplier crossing the unit circle%
\footnote{%
Observe that the proposed function is not suitable to test for a pair of conjugate multipliers crossing the unit circle, as it does not change sign there.
}%
.
We first create the test function as a file \texttt{test\_quadraticre\_bif.m} with the following content.
\begin{lstlisting}
function test = test_quadraticre_bif(gamma, system)
system.par = gamma;
mult = eigTMNpw(system);
% ignore closest to 1 (1 is always present due to linearization)
[~, ind] = min(abs(mu-1));
mult(ind) = [];
test = prod(abs(mult)-1);
end
\end{lstlisting}
We can then find two bifurcations (a Hopf bifurcation at $\gamma = 3.570796208333382$ and a period doubling bifurcation at $\gamma = 4.325285374879225$) with these instructions.
\begin{lstlisting}
fzero(@(gamma) test_quadraticre_bif(gamma, system), 4)
fzero(@(gamma) test_quadraticre_bif(gamma, system), 4.2)
\end{lstlisting}

\subsection{A delay differential equation with a numerical periodic solution computed with DDE-BIFTOOL}
\label{ch:BLV:sec:belairmackey}

We now consider the Bélair--Mackey equation modelling the regulation of mammalian platelet production \cite{BelairMackey1987}
\begin{equation}
\label{ch:BLV:eq:belairmackey}
x'(t) = -\gamma x(t) + q(x(t-\tau_{m})) - q(x(t-\tau_{m}-\tau_{s})) \E^{-\gamma \tau_{s}}
\end{equation}
with $q(x) = q_{0} \theta^{n} x / (\theta^{n} + x^{n})$.
We consider a periodic solution $\bar{x}$ of \eqref{ch:BLV:sec:belairmackey} computed using DDE-BIFTOOL (version 3.2a) with $\gamma = 12$, $q_{0} = 27000$, $n = 2.133$, $\theta = 0.04$, $\tau_{m} = 9$ and $\tau_{s} = 10$; the computed period is $\omega = 18.208035651940627$.
The software repository for \texttt{eigTMNpw} contains the script \texttt{test\_belairmackey\_sol\_db.m} computing the solution and the file \texttt{test\_belairmackey\_sol\_db.mat} containing the result of the computation.
The linearization of \eqref{ch:BLV:eq:belairmackey} around $\bar{x}$ reads
\begin{equation*}
\begin{aligned}
x'(t) &= -\gamma x(t) + q'(\bar{x}(t-\tau_{m})) x(t-\tau_{m}) \\
&\phantom{=}\,\, - \! q'(\bar{x}(t-\tau_{m}-\tau_{s})) \E^{-\gamma \tau_{s}} x(t-\tau_{m}-\tau_{s}).
\end{aligned}
\end{equation*}
As anticipated in Sect.~\ref{ch:BLV:sec:discr}, when using numerical periodic solutions it is essential to follow the corresponding partition of $[0, \omega]$ in discretizing the monodromy operator, especially when this partition is far from being uniform.
This is done in the example below by specifying the parameter \texttt{t} to \texttt{eigTMNpw\_system}.
Figure~\ref{ch:BLV:sec:belairmackey} (left panel) shows the solution $\bar{x}$ and the corresponding partition of $[0, \omega]$.

\begin{figure}[htp]
\centering
\includegraphics[trim=14.113pt 0 0 0]{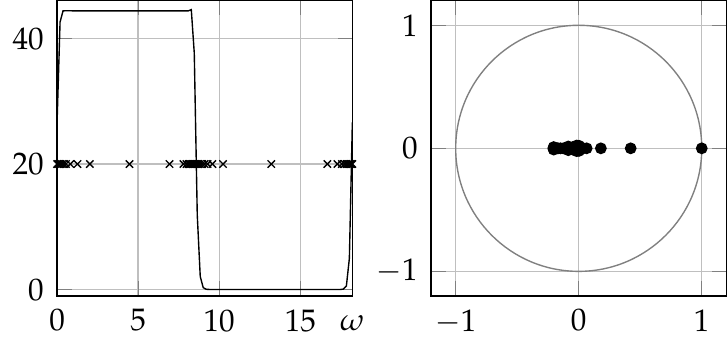}
\caption{%
The solution $\bar{x}$ of \eqref{ch:BLV:eq:belairmackey} with crosses representing the corresponding partition of $[0, \omega]$ (left) and the relevant multipliers (right) computed as described in Sect.~\ref{ch:BLV:sec:belairmackey} with $\gamma = 12$, $q_{0} = 27000$, $n = 2.133$, $\theta = 0.04$, $\tau_{m} = 9$ and $\tau_{s} = 10$ ($\omega = 18.208035651940627$).
}
\label{ch:BLV:fig:belairmackey}
\end{figure}

With the DDE-BIFTOOL point structure \texttt{sol} corresponding to $\bar{x}$ loaded in the workspace, we can construct a structure array \texttt{par} of the parameters, the vector \texttt{mesh} of the endpoints of the partition of $[0, \omega]$ and the piecewise polynomial \texttt{par.sol} representing the numerical solution.
For the latter we can use the helper function \texttt{make\_num\_per\_sol\_db}, available with \texttt{eigTMNpw}.
\begin{lstlisting}
parnames = {'gamma', 'q0', 'n', 'theta', 'tau_m', 'tau_s', 'tau_ms'};
par = cell2struct(num2cell(sol.parameter), parnames, 2);
par.sol = make_num_per_sol_db(sol);
mesh = sol.mesh(1:sol.degree:end);
\end{lstlisting}
We can now define the equation and compute the multipliers (shown in Fig.~\ref{ch:BLV:fig:belairmackey}, right panel) using the following instructions, where for the polynomials in \texttt{eigTMNpw} we use the same degree as the numerical solution.
\begin{lstlisting}
dq = @(x, q0, n, theta) q0 * theta^n ...
  * ((1-n) * x.^n + theta^n) ./ (theta^n + x.^n)^2;
system = eigTMNpw_system(...
  [0, 1], ...                  % dimensions [dX, dY]
  [par.tau_m, par.tau_ms], ... % delays
  sol.period, ...              % period (or length of time evolution)
  0, ...                       % initial time
  'AYY', @(t, par) -par.gamma, ...
  'BYY', {...
    @(t, par) dq(par.sol(1, t-par.tau_m), par.q0, par.n, par.theta), ...
    @(t, par) dq(par.sol(1, t-par.tau_ms), par.q0, par.n, par.theta) ...
              * exp(-par.gamma*par.tau_s)}, ...
  'par', par, ...
  't', mesh);
method = eigTMNpw_method(sol.degree);
mult = eigTMNpw(system, method);
\end{lstlisting}
The computed dominant multipliers are $0.999974005170910$, $0.422325480377944$ and $-0.204620549659091\pm0.004612509701759\I$.

\subsection{A coupled equation with a numerical periodic solution computed with MatCont}
\label{ch:BLV:sec:logisticdaphniadelayed}

For the next example we consider a simplified version of the \emph{Daphnia} model \cite{DiekmannGyllenbergMetzNakaokaDeRoos2010} with explicit terms for the survival probability, a fixed maturation age and a consumer-free resource dynamics of delayed logistic type, namely the nonlinear coupled equation
\begin{equation}
\label{ch:BLV:eq:logisticdaphniadelayed}
\left\{
\begin{aligned}
& b(t) = \beta S(t) \int_{\bar{a}}^{\tau} b(t-a) \D a, \\
& S'(t) = r S(t) \left(1 - \frac{S(t-\tau)}{K}\right) - \gamma S(t) \int_{\bar{a}}^{\tau} b(t-a) \D a,
\end{aligned}
\right.
\end{equation}
where all parameters are positive and $\bar{a} < \tau$.
We consider a periodic solution $(\bar{b}, \bar{S})$ of \eqref{ch:BLV:eq:logisticdaphniadelayed} computed using MatCont%
\footnote{\url{https://matcont.sourceforge.io/}, version 7.2.}
\cite{DhoogeGovaertsKuznetsovMeijerSautois2008} discretizing the equation as a system of ordinary differential equations according to \cite{BredaDiekmannGyllenbergScarabelVermiglio2016}, with $\beta = 2$, $\bar{a} = 3$, $r = 0.3$, $K = 1$, $\gamma = 1$ and $\tau = 4$; the computed period is $\omega = 23.133253862004800$.
The software repository for \texttt{eigTMNpw} contains the script \texttt{test\_logisticdaphniadelayed\_sol\_mc.m} computing the solution and the file \texttt{test\_logisticdaphniadelayed\_sol\_db.mat} containing the result of the computation.
Figure~\ref{ch:BLV:sec:logisticdaphniadelayed} (left panel) shows the solution $(\bar{b}, \bar{S})$ and the corresponding partition of $[0, \omega]$.
To study the stability of $(\bar{b}, \bar{S})$ we consider the linear coupled equation
\begin{equation*}
\left\{
\begin{aligned}
& b(t) = \beta \bar{S}(t) \int_{\bar{a}}^{\tau} b(t-a) \D a + \beta \int_{\bar{a}}^{\tau} \bar{b}(t-a) \D a \, S(t), \\
& S'(t) = \left[r \left(1 - \frac{\bar{S}(t-\tau)}{K}\right) - \gamma \int_{\bar{a}}^{\tau} \bar{b}(t-a) \D a\right] S(t) \\
& \phantom{S'(t) ={}} - \frac{r \bar{S}(t)}{K} S(t-\tau) - \gamma \bar{S}(t) \int_{\bar{a}}^{\tau} b(t-a) \D a.
\end{aligned}
\right.
\end{equation*}

\begin{figure}[htp]
\centering
\includegraphics[trim=22.442pt 0 0 0]{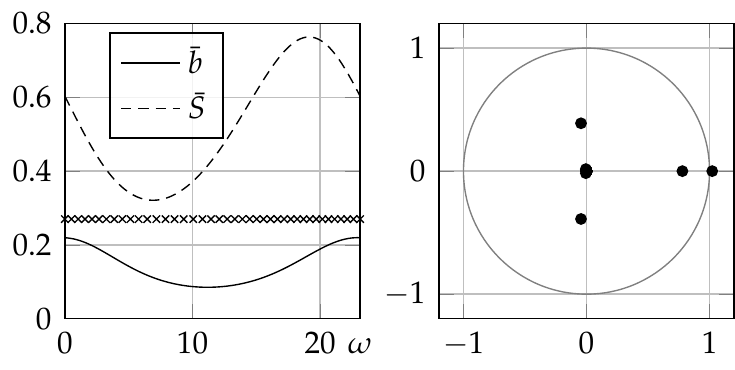}
\caption{%
The solution $(\bar{b}, \bar{S})$ of \eqref{ch:BLV:eq:logisticdaphniadelayed} with crosses representing the corresponding partition of $[0, \omega]$ (left) and the relevant multipliers (right) computed as described in Sect.~\ref{ch:BLV:sec:logisticdaphniadelayed} with $\beta = 2$, $\bar{a} = 3$, $r = 0.3$, $K = 1$, $\gamma = 1$ and $\tau = 4$ ($\omega = 23.133253862004800$).
}
\label{ch:BLV:fig:logisticdaphniadelayed}
\end{figure}

With MatCont's solution loaded in the workspace, we can recover the \texttt{period}, and construct a structure array \texttt{par} of the parameters, the vector \texttt{mesh} of the endpoints of the partition of $[0, \omega]$ and the piecewise polynomial \texttt{par.sol} representing the numerical solution.
For the latter we can use the helper function \texttt{make\_num\_per\_sol\_mc}, available with \texttt{eigTMNpw}.
\begin{lstlisting}
parnames = {'beta', 'abar', 'r', 'K', 'gamma', 'aux', 'tau', 'M'};
par = [parnames; num2cell(slc.data.parametervalues')];
par = struct(par{:});
par.sol = make_num_per_sol_mc(dX, dY, xlc, slc, flc);
mesh = flc(1:slc.data.ntst+1);
period = xlc(end-1);
\end{lstlisting}
We can now define the equation and compute the multipliers (shown in Fig.~\ref{ch:BLV:fig:logisticdaphniadelayed}, right panel) using the following instructions, where for the polynomials in \texttt{eigTMNpw} we use the same degree as the numerical solution.
\begin{lstlisting}
system = eigTMNpw_system(...
  [dX, dY], ...            % dimensions [dX, dY]
  [par.abar, par.tau], ... % delays
  period, ...              % period (or length of time evolution)
  0, ...                   % initial time
  'AXY', @(t, par) par.beta ...
         * integral(@(s) par.sol(1, t-s), par.abar, par.tau), ...
  'AYY', @(t, par) par.r * (1 - par.sol(2, t-par.tau) / par.K) ...
         - par.gamma ...
         * integral(@(s) par.sol(1, t-s), par.abar, par.tau), ...
  'BYY', {[], @(t, par) -par.r / par.K * par.sol(2, t)}, ...
  'CXX', {[], @(t, theta, par) par.beta * par.sol(2, t)}, ...
  'CYX', {[], @(t, theta, par) -par.gamma * par.sol(2, t)}, ...
  'par', par, ...
  't', mesh);
method = eigTMNpw_method(slc.data.ncol);
mult = eigTMNpw(system, method);
\end{lstlisting}
The computed dominant multipliers are $1.021824635351366$, $0.779435823328075$ and $-0.044734353116124\pm0.388939410834478\I$.

\subsection{Computing a stability chart}
\label{ch:BLV:sec:hayes}

\begin{figure}[htp]
\centering
\includegraphics[trim=19.418pt 0 0 0]{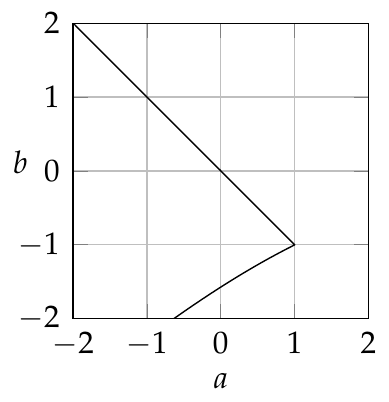}
\caption{%
Stability chart for \eqref{ch:BLV:eq:hayes} with $\tau = 1$ as computed in Sect.~\ref{ch:BLV:sec:hayes}.
The null equilibrium is asymptotically stable to the left of the curve.
}
\label{ch:BLV:fig:hayes}
\end{figure}

As a last example we want to compute a stability chart for the Hayes equation
\begin{equation}
\label{ch:BLV:eq:hayes}
x'(t) = a x(t) + b x(t-\tau)
\end{equation}
with $a, b \in [-2, 2]$ and $\tau = 1$.
In order to check whether all multipliers are inside the unit circle (asymptotically stable null equilibrium) or any are outside (unstable null equilibrium), we compute the difference between the maximum magnitude of the multipliers and $1$ and find the curve of level $0$ of the resulting function.
We first create the test function as a file \texttt{test\_hayes\_stab.m} with the following content.
\begin{lstlisting}
function test = test_hayes_stab(a, b, system)
system.par = [a, b];
mult = eigTMNpw(system);
test = max(abs(mult)) - 1;
end
\end{lstlisting}
Then we define the equation with the following instruction.
\begin{lstlisting}
system = eigTMNpw_system(...
  [0, 1], ... % dimensions [dX, dY]
  1, ...      % delays
  1, ...      % period (or length of time evolution)
  0, ...      % initial time
  'AYY', @(t, par) par(1), ...
  'BYY', {@(t, par) par(2)});
\end{lstlisting}
Finally, we can find the curve of level $0$%
\footnote{%
Observe that we actually compute the curve of level $10^{-3}$, since the curve of level $0$ poses some numerical problems, probably due to simmetry; however, the curves of levels $10^{-3}$ and $-10^{-3}$ are indistinguishable for the purpose of this illustration.
}%
,
shown in Fig.~\ref{ch:BLV:fig:hayes}, using LEVEL%
\footnote{%
\url{http://cdlab.uniud.it/software}
}
\cite{BredaMasetVermiglio2009b}.
\begin{lstlisting}
level(@(a, b) test_hayes_stab(a, b, system), ...
  1e-3, ...           % level
  [-2, 2, -2, 2], ... % domain
  1e-1);              % tolerance
\end{lstlisting}

\section{Advanced topics}
\label{ch:BLV:sec:adv}

We complete the work with the following comments on two additional use cases, namely deriving a formulation as a generalized eigenvalue problem (Sect.~\ref{ch:BLV:sec:gep}) and defining new families of collocation nodes and new quadrature formulas (Sect.~\ref{ch:BLV:sec:new}).

\subsection{Generalized eigenvalue problems}
\label{ch:BLV:sec:gep}

As mentioned in Sect.~\ref{ch:BLV:sec:method}, \texttt{eigTMNpw} outputs in order the vector of the eigenvalues of $T^{\bullet}$, the matrix $T^{\bullet}$ and four other matrices $T^{\bullet}_{1}$, $T^{\bullet}_{2}$, $U^{\bullet}_{1}$ and $U^{\bullet}_{2}$.
The latter are defined as follows.

Recalling \eqref{ch:BLV:eq:V}, let
\begin{alignat*}{2}
V^{-}&\colon X \times Y \to X^{\pm} \times Y^{\pm},
&\quad
V^{+}&\colon X^{+} \times Y^{+} \to X^{\pm} \times Y^{\pm},
\\
T^{\bullet}_{1} &\colon X^{\bullet} \times Y^{\bullet} \to X^{\bullet} \times Y^{\bullet},
&\quad
T^{\bullet}_{2} &\colon X^{+\bullet} \times Y^{+\bullet} \to X^{\bullet} \times Y^{\bullet},
\\
U^{\bullet}_{1} &\colon X^{\bullet} \times Y^{\bullet} \to X^{+\bullet} \times Y^{+\bullet},
&\quad
U^{\bullet}_{2} &\colon X^{+\bullet} \times Y^{+\bullet} \to X^{+\bullet} \times Y^{+\bullet},
\end{alignat*}
be defined as
\begin{alignat*}{2}
V^{-}(\phi, \psi) &\coloneqq V((\phi, \psi), (0_{X^{+}}, 0_{Y^{+}})),
&\quad
V^{+}(w, z) &\coloneqq V((0_{X}, 0_{Y}), (w, z)),
\\
T^{\bullet}_{1} (\Phi, \Psi) &\coloneqq R (V^{-} P (\Phi, \Psi))_{\omega},
&\quad
T^{\bullet}_{2} (W, Z) &\coloneqq R (V^{+} P^{+} (W, Z))_{\omega},
\\
U^{\bullet}_{1} (\Phi, \Psi) &\coloneqq R^{+} \mathcal{F}_{s} V^{-} P (\Phi, \Psi),
&\quad
U^{\bullet}_{2} (W, Z) &\coloneqq R^{+} \mathcal{F}_{s} V^{+} P^{+} (W, Z).
\end{alignat*}
Recalling \eqref{ch:BLV:eq:discrete-T-as-V}--\eqref{ch:BLV:eq:discrete-fixed-point}, observing that $V((\phi, \psi), (w,z)) = V^{-}(\phi, \psi) + V^{+}(w, z)$ and thanks to the linearity of $\mathcal{F}_{s}$, the operator $T^{\bullet}$ can be rewritten as
\begin{equation*}
T^{\bullet} = T^{\bullet}_{1} + T^{\bullet}_{2} (I_{X^{+\bullet} \times Y^{+\bullet}} - U^{\bullet}_{2})^{-1} U^{\bullet}_{1}
\end{equation*}
where the invertibility of $I_{X^{+\bullet} \times Y^{+\bullet}} - U^{\bullet}_{2}$ is ensured under reasonable hypotheses, see \cite[Sect.~4.2]{BredaLiessi2020}.

This reformulation simplifies the construction of the matrix representation of $T^{\bullet}$, needed for the implementation of the method.
In our formulation we explicitly compute the matrix, allowing us to approximate any evolution operator, whose eigenvalues can be determined, if needed, by solving a standard eigenvalue problem $T^{\bullet} (\Phi, \Psi) = \mu (\Phi, \Psi)$.
However, using the matrices $T^{\bullet}_{1}$, $T^{\bullet}_{2}$, $U^{\bullet}_{1}$ and $U^{\bullet}_{2}$ separately, we can reformulate the problem as the generalized eigenvalue problem
\small
\begin{equation*}
\left[
\begin{pmatrix}
T^{\bullet}_{1} & T^{\bullet}_{2} \\
U^{\bullet}_{1} & U^{\bullet}_{2}
\end{pmatrix}
- I_{(X^{\bullet} \times Y^{\bullet}) \times (X^{+\bullet} \times Y^{+\bullet})}
\right]
\begin{pmatrix}
(\Phi, \Psi) \\
(W, Z)
\end{pmatrix}
=
(\mu - 1)
\begin{pmatrix}
I_{X^{\bullet} \times Y^{\bullet}} & 0_{\bullet\,\leftarrow\,+\bullet} \\
0_{+\bullet\,\leftarrow\,\bullet} & 0_{X^{+\bullet} \times Y^{+\bullet}}
\end{pmatrix}
\begin{pmatrix}
(\Phi, \Psi) \\
(W, Z)
\end{pmatrix}
,
\end{equation*}
\normalsize
with the $0$ operators in the antidiagonal defined as $0_{\bullet\,\leftarrow\,+\bullet} \colon X^{+\bullet} \times Y^{+\bullet} \to X^{\bullet} \times Y^{\bullet}$ and $0_{+\bullet\,\leftarrow\,\bullet} \colon X^{\bullet} \times Y^{\bullet} \to X^{+\bullet} \times Y^{+\bullet}$ (the resulting generalized eigenvalue problem is similar to the one proposed in \cite{BorgioliHajduInspergerStepanMichiels2020}).
Since in the piecewise approach the matrices $T^{\bullet}_{1}$, $T^{\bullet}_{2}$ are sparse and the matrices $U^{\bullet}_{1}$ and $U^{\bullet}_{2}$ often are as well, the generalized formulation allows us to possibly exploit their structure for efficiency.
Our tests indicate that the generalized approach becomes computationally convenient for large dimensions of the problem.

We can obtain the four matrices, along with the multipliers computed via the standard eigenvalue problem and the full matrix $T^{\bullet}$, with the following instruction.
\begin{lstlisting}
[mult, T, T1, T2, U1, U2] = eigTMNpw(system, ...);
\end{lstlisting}
Note that our code constructs the matrices in full form, hence taking full advantage of the structure would require to reimplement the discretization method with an eye to sparsity.

\subsection{User defined nodes and quadrature}
\label{ch:BLV:sec:new}

As a last comment, users might want to use a family of collocation nodes and a quadrature method of their choice.
Defining new options for them is quite simple, requiring changes only in \texttt{eigTMNpw\_method} to add new keywords and corresponding cases in \texttt{validCollocationFamilies} and \texttt{switch method.CollocationFamily}, and in \texttt{validQuadratureMethods} and \texttt{switch method.QuadratureMethod}, respectively.
For a collocation family the nodes themselves (given in ascending order in $[0, 1]$) and the corresponding weights for the barycentric Lagrange interpolation must be defined, while for a quadrature method a function taking the function handle of the integrand and the integration endpoints must be provided, possibly supporting options via the \texttt{QuadratureParameters} named parameter of \texttt{eigTMNpw\_method}.
For example, Chebyshev extrema are defined as follows in \texttt{eigTMNpw\_method}, with \texttt{method.c} and \texttt{method.bw} being the required elements.
\begin{lstlisting}
switch method.CollocationFamily
  case 'cheb2'
    % Chebyshev type II (extrema) points in [0, 1] (in ascending order)
    method.c = (1 - cos((0:method.M)*pi/method.M)) / 2;
    % and corresponding barycentric weights
    method.bw = [1/2, ones(1, method.M-1), 1/2] .* (-1) .^ (0:method.M);
\end{lstlisting}

\section*{Acknowledgements}

The authors are members of INdAM Research group GNCS and of UMI Research group ``Mo\-del\-li\-sti\-ca socio-epidemiologica''.
This work was partially supported by the Italian Ministry of University and Research (MUR) through the PRIN 2020 project (No.\ 2020JLWP23) ``Integrated Mathematical Approaches to Socio-Epidemiological Dynamics'' (CUP: E15F21005420006).
The work of Davide Liessi was partially supported by Finanziamento Giovani Ricercatori 2020--2021 of INdAM Research group GNCS.

{\sloppy
\printbibliography

@ARTICLE{Ando2021,
  AUTHOR = {Andò, Alessia},
  DATE = {2021-09-08},
  DOI = {10.1002/cmm4.1190},
  JOURNALTITLE = {Computational and Mathematical Methods},
  PAGES = {e1190},
  SHORTJOURNAL = {Comput. Math. Methods},
  TITLE = {Convergence of collocation methods for solving periodic boundary value problems for renewal equations defined through finite-dimensional boundary conditions},
}

@ARTICLE{AndoBreda2020b,
  AUTHOR = {Andò, Alessia and Breda, Dimitri},
  DATE = {2020-10-27},
  DOI = {10.1137/19M1295015},
  JOURNALTITLE = {SIAM Journal on Numerical Analysis},
  NOTE = {full-length version at \ifhyperref{\href{https://arxiv.org/abs/2008.07604}{\nolinkurl{arxiv:2008.07604}}}{\nolinkurl{arxiv:2008.07604}} [math.NA]},
  NUMBER = {5},
  PAGES = {3010--3039},
  SHORTJOURNAL = {SIAM J. Numer. Anal.},
  TITLE = {Convergence analysis of collocation methods for computing periodic solutions of retarded functional differential equations},
  VOLUME = {58},
}

@ARTICLE{AndoBreda,
  AUTHOR = {Andò, Alessia and Breda, Dimitri},
  DATE = {2021},
  PUBSTATE = {submitted},
  TITLE = {Numerical computation of periodic solutions of renewal equations from structured population dynamics},
}

@INPROCEEDINGS{ArinoVanDenDriessche2006,
  AUTHOR = {Arino, O. and van den Driessche, P.},
  EDITOR = {Arino, O. and Hbid, M. L. and Ait Dads, E.},
  LOCATION = {Dordrecht},
  PUBLISHER = {Springer},
  BOOKSUBTITLE = {Proceedings of the NATO Advanced Study Institute},
  BOOKTITLE = {Delay Differential Equations and Applications},
  DATE = {2006},
  DOI = {10.1007/1-4020-3647-7},
  EVENTDATE = {2002-09},
  NUMBER = {205},
  SERIES = {NATO Science Series II: Mathematics, Physics and Chemistry},
  SHORTSERIES = {NATO Sci. Ser. II Math. Phys. Chem.},
  SUBTITLE = {Modeling and Numerical Considerations},
  TITLE = {Time Delays in Epidemic Models},
  VENUE = {Marrakech, Morocco},
}

@INCOLLECTION{Bader1985,
  AUTHOR = {Bader, G.},
  EDITOR = {Ascher, U. M. and Russell, R. D.},
  LOCATION = {Boston},
  PUBLISHER = {Birkhäuser},
  BOOKTITLE = {Numerical Boundary Value ODEs},
  DATE = {1985},
  DOI = {10.1007/978-1-4612-5160-6_13},
  NUMBER = {5},
  PAGES = {227--243},
  SERIES = {Progress in Scientific Computing},
  SHORTSERIES = {Progr. Sci. Comput.},
  TITLE = {Solving boundary value problems for functional differential equations by collation},
}

@ARTICLE{BelairMackey1987,
  AUTHOR = {Bélair, Jacques and Mackey, Michael C.},
  DATE = {1987-07},
  DOI = {10.1111/j.1749-6632.1987.tb48740.x},
  JOURNALTITLE = {Annals of the New York Academy of Sciences},
  NUMBER = {1},
  PAGES = {280--282},
  SHORTJOURNAL = {Ann. New York Acad. Sci.},
  TITLE = {A model for the regulation of mammalian platelet production},
  VOLUME = {504},
}

@ARTICLE{BorgioliHajduInspergerStepanMichiels2020,
  AUTHOR = {Borgioli, Francesco and Hajdu, David and Insperger, Tamás and Stépán, Gábor and Michiels, Wim},
  DATE = {2020-08-31},
  DOI = {10.1002/nme.6368},
  JOURNALTITLE = {International Journal for Numerical Methods in Engineering},
  NUMBER = {16},
  PAGES = {3505--3528},
  SHORTJOURNAL = {Internat. J. Numer. Methods Engrg.},
  TITLE = {Pseudospectral method for assessing stability robustness for linear time-periodic delayed dynamical systems},
  VOLUME = {121},
}

@ARTICLE{BredaDiekmannDeGraafPuglieseVermiglio2012,
  AUTHOR = {Breda, Dimitri and Diekmann, Odo and de Graaf, W. F. and Pugliese, A . and Vermiglio, Rossana},
  DATE = {2012},
  DOI = {10.1080/17513758.2012.716454},
  JOURNALTITLE = {Journal of Biological Dynamics},
  NUMBER = {sup2},
  PAGES = {103--117},
  SHORTJOURNAL = {J. Biol. Dyn.},
  TITLE = {On the formulation of epidemic models (an appraisal of {K}ermack and {M}c{K}endrick)},
  VOLUME = {6},
}

@ARTICLE{BredaDiekmannGyllenbergScarabelVermiglio2016,
  AUTHOR = {Breda, Dimitri and Diekmann, Odo and Gyllenberg, Mats and Scarabel, Francesca and Vermiglio, Rossana},
  DATE = {2016},
  DOI = {10.1137/15M1040931},
  JOURNALTITLE = {SIAM Journal on Applied Dynamical Systems},
  NUMBER = {1},
  PAGES = {1--23},
  SHORTJOURNAL = {SIAM J. Appl. Dyn. Syst.},
  TITLE = {Pseudospectral discretization of nonlinear delay equations: new prospects for numerical bifurcation analysis},
  VOLUME = {15},
}

@ARTICLE{BredaDiekmannLiessiScarabel2016,
  AUTHOR = {Breda, Dimitri and Diekmann, Odo and Liessi, Davide and Scarabel, Francesca},
  DATE = {2016-09-12},
  DOI = {10.14232/ejqtde.2016.1.65},
  EID = {65},
  JOURNALTITLE = {Electronic Journal of Qualitative Theory of Differential Equations},
  PAGES = {1--24},
  SHORTJOURNAL = {Electron. J. Qual. Theory Differ. Equ.},
  TITLE = {Numerical bifurcation analysis of a class of nonlinear renewal equations},
  VOLUME = {2016},
}

@ARTICLE{BredaLiessi2018,
  AUTHOR = {Breda, Dimitri and Liessi, Davide},
  DATE = {2018},
  DOI = {10.1137/17M1140534},
  JOURNALTITLE = {SIAM Journal on Numerical Analysis},
  NUMBER = {3},
  PAGES = {1456--1481},
  SHORTJOURNAL = {SIAM J. Numer. Anal.},
  TITLE = {Approximation of eigenvalues of evolution operators for linear renewal equations},
  VOLUME = {56},
}

@ARTICLE{BredaLiessi2020,
  AUTHOR = {Breda, Dimitri and Liessi, Davide},
  DATE = {2020-11},
  DOI = {10.1007/s11587-020-00513-9},
  IDS = {BredaLiessi2020b},
  JOURNALTITLE = {Ricerche di Matematica},
  NUMBER = {2},
  PAGES = {457--481},
  SHORTJOURNAL = {Ric. Mat.},
  TITLE = {Approximation of eigenvalues of evolution operators for linear coupled renewal and retarded functional differential equations},
  VOLUME = {69},
}

@ARTICLE{BredaLiessi2021,
  AUTHOR = {Breda, Dimitri and Liessi, Davide},
  DATE = {2021-06},
  DOI = {10.1007/s10884-020-09826-7},
  IDS = {BredaLiessi2020a},
  JOURNALTITLE = {Journal of Dynamics and Differential Equations},
  NUMBER = {2},
  PAGES = {457--481},
  SHORTJOURNAL = {J. Dynam. Differential Equations},
  TITLE = {Floquet theory and stability of periodic solutions of renewal equations},
  VOLUME = {33},
}

@ARTICLE{BredaLiessiVermiglio,
  AUTHOR = {Breda, Dimitri and Liessi, Davide and Vermiglio, Rossana},
  DATE = {2022},
  JOURNALTITLE = {Journal of Computational Dynamics},
  PUBSTATE = {forthcoming},
  SHORTJOURNAL = {J. Comput. Dyn.},
  TITLE = {Piecewise discretization of monodromy operators of delay equations on adapted meshes},
}

@ARTICLE{BredaMasetVermiglio2009b,
  AUTHOR = {Breda, Dimitri and Maset, Stefano and Vermiglio, Rossana},
  DATE = {2009-12},
  DOI = {10.1007/s11075-009-9303-2},
  JOURNALTITLE = {Numerical Algorithms},
  NUMBER = {4},
  PAGES = {605--628},
  SHORTJOURNAL = {Numer. Algorithms},
  TITLE = {An adaptive algorithm for efficient computation of level curves of surfaces},
  VOLUME = {52},
}

@ARTICLE{BredaMasetVermiglio2012,
  AUTHOR = {Breda, Dimitri and Maset, Stefano and Vermiglio, Rossana},
  DATE = {2012},
  DOI = {10.1137/100815505},
  JOURNALTITLE = {SIAM Journal on Numerical Analysis},
  NUMBER = {3},
  PAGES = {1456--1483},
  SHORTJOURNAL = {SIAM J. Numer. Anal.},
  TITLE = {Approximation of eigenvalues of evolution operators for linear retarded functional differential equations},
  VOLUME = {50},
}

@BOOK{BredaMasetVermiglio2015,
  AUTHOR = {Breda, Dimitri and Maset, Stefano and Vermiglio, Rossana},
  LOCATION = {New York},
  PUBLISHER = {Springer},
  DATE = {2015},
  DOI = {10.1007/978-1-4939-2107-2},
  SERIES = {SpringerBriefs in Control, Automation and Robotics},
  SHORTSERIES = {SpringerBriefs Control Autom. Robot.},
  SUBTITLE = {A Numerical Approach with MATLAB},
  TITLE = {Stability of Linear Delay Differential Equations},
}

@ARTICLE{Bueler2007,
  AUTHOR = {Bueler, Ed},
  DATE = {2007},
  DOI = {10.1137/050633330},
  JOURNALTITLE = {SIAM Journal on Numerical Analysis},
  NUMBER = {6},
  PAGES = {2510--2536},
  SHORTJOURNAL = {SIAM J. Numer. Anal.},
  TITLE = {Error bounds for approximate eigenvalues of periodic-coefficient linear delay differential equations},
  VOLUME = {45},
}

@ARTICLE{ButcherBobrenkov2011,
  AUTHOR = {Butcher, Eric A. and Bobrenkov, Oleg A.},
  DATE = {2011-03},
  DOI = {10.1016/j.cnsns.2010.05.037},
  JOURNALTITLE = {Communications in Nonlinear Science and Numerical Simulation},
  NUMBER = {3},
  PAGES = {1541--1554},
  SHORTJOURNAL = {Commun. Nonlinear Sci. Numer. Simul.},
  TITLE = {On the {C}hebyshev spectral continuous time approximation for constant and periodic delay differential equations},
  VOLUME = {16},
}

@ARTICLE{ButcherMaBuelerAverinaSzabo2004,
  AUTHOR = {Butcher, Eric A. and Ma, Haitao and Bueler, Ed and Averina, Victoria and Szabó, Zsolt},
  DATE = {2004-02-21},
  DOI = {10.1002/nme.894},
  JOURNALTITLE = {International Journal for Numerical Methods in Engineering},
  NUMBER = {7},
  PAGES = {895--922},
  SHORTJOURNAL = {Internat. J. Numer. Methods Engrg.},
  TITLE = {Stability of linear time-periodic delay-differential equations via {C}hebyshev polynomials},
  VOLUME = {59},
}

@ARTICLE{ClenshawCurtis1960,
  AUTHOR = {Clenshaw, C. W. and Curtis, A. R.},
  DATE = {1960-12},
  DOI = {10.1007/BF01386223},
  JOURNALTITLE = {Numerische Mathematik},
  PAGES = {197--205},
  SHORTJOURNAL = {SIAM J. Numer. Anal.},
  TITLE = {A method for numerical integration on an automatic computer},
  VOLUME = {2},
}

@ARTICLE{DhoogeGovaertsKuznetsovMeijerSautois2008,
  AUTHOR = {Dhooge, A. and Govaerts, W. and Kuznetsov, Yuri A. and Meijer, H. G. E. and Sautois, B.},
  DATE = {2008-04},
  DOI = {10.1080/13873950701742754},
  JOURNALTITLE = {Mathematical and Computer Modelling of Dynamical Systems},
  NUMBER = {2},
  PAGES = {147--175},
  SHORTJOURNAL = {Math. Comput. Model. Dyn. Syst.},
  TITLE = {New features of the software {MatCont} for bifurcation analysis of dynamical systems},
  VOLUME = {14},
}

@ARTICLE{DiekmannGettoGyllenberg2008,
  AUTHOR = {Diekmann, Odo and Getto, Philipp and Gyllenberg, Mats},
  DATE = {2008},
  DOI = {10.1137/060659211},
  JOURNALTITLE = {SIAM Journal on Mathematical Analysis},
  NUMBER = {4},
  PAGES = {1023--1069},
  SHORTJOURNAL = {SIAM J. Math. Anal.},
  TITLE = {Stability and bifurcation analysis of {V}olterra functional equations in the light of suns and stars},
  VOLUME = {39},
}

@ARTICLE{DiekmannGyllenbergMetzNakaokaDeRoos2010,
  AUTHOR = {Diekmann, Odo and Gyllenberg, Mats and Metz, J. A. J. and Nakaoka, Shinji and de Roos, André Marc},
  DATE = {2010-08},
  DOI = {10.1007/s00285-009-0299-y},
  JOURNALTITLE = {Journal of Mathematical Biology},
  NUMBER = {2},
  PAGES = {277--318},
  SHORTJOURNAL = {J. Math. Biol.},
  TITLE = {\emph{Daphnia} revisited: local stability and bifurcation theory for physiologically structured population models explained by way of an example},
  VOLUME = {61},
}

@BOOK{DiekmannVanGilsVerduynLunelWalther1995,
  AUTHOR = {Diekmann, Odo and van Gils, Stephan A. and Verduyn Lunel, Sjoerd M. and Walther, Hans-Otto},
  LOCATION = {New York},
  PUBLISHER = {Springer},
  DATE = {1995},
  DOI = {10.1007/978-1-4612-4206-2},
  NUMBER = {110},
  SERIES = {Applied Mathematical Sciences},
  SHORTSERIES = {Appl. Math. Sci.},
  SUBTITLE = {Functional-, Complex- and Nonlinear Analysis},
  TITLE = {Delay Equations},
}

@ARTICLE{EngelborghsLuzyaninaIntHoutRoose2001,
  AUTHOR = {Engelborghs, Koen and Luzyanina, Tatyana and in 't Hout, Karel J. and Roose, Dirk},
  DATE = {2001},
  DOI = {10.1137/S1064827599363381},
  JOURNALTITLE = {SIAM Journal on Scientific Computing},
  NUMBER = {5},
  PAGES = {1593--1609},
  SHORTJOURNAL = {SIAM J. Sci. Comput.},
  TITLE = {Collocation methods for the computation of periodic solutions of delay differential equations},
  VOLUME = {22},
}

@ARTICLE{EngelborghsLuzyaninaRoose2002,
  AUTHOR = {Engelborghs, Koen and Luzyanina, Tatyana and Roose, Dirk},
  DATE = {2002},
  DOI = {10.1145/513001.513002},
  JOURNALTITLE = {ACM Transactions on Mathematical Software},
  NUMBER = {1},
  PAGES = {1--21},
  SHORTJOURNAL = {ACM Trans. Math. Softw.},
  TITLE = {Numerical Bifurcation Analysis of Delay Differential Equations Using {DDE-BIFTOOL}},
  VOLUME = {28},
}

@BOOK{Erneux2009,
  AUTHOR = {Erneux, Thomas},
  LOCATION = {New York},
  PUBLISHER = {Springer},
  DATE = {2009},
  DOI = {10.1007/978-0-387-74372-1},
  NUMBER = {3},
  SERIES = {Surveys and Tutorials in the Applied Mathematical Sciences},
  SHORTSERIES = {Surveys Tutor. Appl. Math. Sci.},
  TITLE = {Applied Delay Differential Equations},
}

@BOOK{Fridman2014,
  AUTHOR = {Fridman, Emilia},
  LOCATION = {Basel},
  PUBLISHER = {Birkhäuser},
  DATE = {2014},
  DOI = {10.1007/978-3-319-09393-2},
  SERIES = {Systems \& Control: Foundations \& Applications},
  SHORTSERIES = {Systems Control Found. Appl.},
  TITLE = {Introduction to Time-Delay Systems: Analysis and Control},
}

@BOOK{GuKharitonovChen2003,
  AUTHOR = {Gu, Keqin and Kharitonov, Vladimir L. and Chen, Jie},
  LOCATION = {Boston, MA},
  PUBLISHER = {Birkhäuser},
  DATE = {2003},
  DOI = {10.1007/978-1-4612-0039-0},
  SERIES = {Control Engineering},
  SHORTSERIES = {Control Eng.},
  TITLE = {Stability of Time-Delay Systems},
}

@ARTICLE{Hayes1950,
  AUTHOR = {Hayes, N. D.},
  DATE = {1950-07},
  DOI = {10.1112/jlms/s1-25.3.226},
  JOURNALTITLE = {Journal of the London Mathematical Society},
  NUMBER = {3},
  PAGES = {226--232},
  SHORTJOURNAL = {J. Lond. Math. Soc.},
  TITLE = {Roots of the transcendental equation associated with a certain difference-differential equation},
  VOLUME = {s1-25},
}

@BOOK{InspergerStepan2011,
  AUTHOR = {Insperger, Tamás and Stépán, Gábor},
  LOCATION = {New York},
  PUBLISHER = {Springer},
  DATE = {2011},
  DOI = {10.1007/978-1-4614-0335-7},
  NUMBER = {178},
  SERIES = {Applied Mathematical Sciences},
  SHORTSERIES = {Appl. Math. Sci.},
  SUBTITLE = {Stability and Engineering Applications},
  TITLE = {Semi-Discretization for Time-Delay Systems},
}

@THESIS{Jarlebring2008,
  AUTHOR = {Jarlebring, Elias},
  INSTITUTION = {Institut Computational Mathematics, Carl-Friedrich-Gauß-Fakultät, Technische Universität Braunschweig},
  DATE = {2008-05-28},
  TITLE = {The spectrum of delay-differential equations: numerical methods, stability and perturbation},
  TYPE = {phdthesis},
}

@BOOK{KolmanovskiiMyshkis1999,
  AUTHOR = {Kolmanovskii, V. and Myshkis, A.},
  LOCATION = {Dordrecht},
  PUBLISHER = {Kluwer Academic Publishers},
  DATE = {1999},
  DOI = {10.1007/978-94-017-1965-0},
  NUMBER = {463},
  SERIES = {Mathematics and Its Applications},
  SHORTSERIES = {Math. Appl.},
  TITLE = {Introduction to the Theory and Applications of Functional Differential Equations},
}

@BOOK{Kuang1993,
  AUTHOR = {Kuang, Yang},
  LOCATION = {San Diego},
  PUBLISHER = {Academic Press},
  DATE = {1993},
  NUMBER = {191},
  SERIES = {Mathematics in Science and Engineering},
  SHORTSERIES = {Math. Sci. Eng.},
  SUBTITLE = {With Applications in Population Dynamics},
  TITLE = {Delay Differential Equations},
}

@ARTICLE{LehotzkyInsperger2016,
  AUTHOR = {Lehotzky, Dávid and Insperger, Tamás},
  DATE = {2016-11-09},
  DOI = {10.1002/nme.5225},
  JOURNALTITLE = {International Journal of Numerical Methods in Engineering},
  NUMBER = {6},
  PAGES = {588--613},
  SHORTJOURNAL = {Internat. J. Numer. Methods Engrg.},
  TITLE = {A pseudospectral tau approximation for time delay systems and its comparison with other weighted-residual-type methods},
  VOLUME = {108},
}

@BOOK{MacDonald1978,
  AUTHOR = {MacDonald, Norman},
  LOCATION = {Berlin, Heidelberg},
  PUBLISHER = {Springer},
  DATE = {1978},
  DOI = {10.1007/978-3-642-93107-9},
  NUMBER = {27},
  SERIES = {Lecture Notes in Biomathematics},
  SHORTSERIES = {Lect. Notes Biomath.},
  TITLE = {Time Lags in Biological Models},
}

@COLLECTION{MetzDiekmann1986,
  EDITOR = {Metz, J. A. J. and Diekmann, Odo},
  LOCATION = {Berlin, Heidelberg},
  PUBLISHER = {Springer},
  DATE = {1986},
  DOI = {10.1007/978-3-662-13159-6},
  NUMBER = {68},
  SERIES = {Lecture Notes in Biomathematics},
  SHORTSERIES = {Lect. Notes Biomath.},
  TITLE = {The Dynamics of Physiologically Structured Populations},
}

@BOOK{MichielsNiculescu2014,
  AUTHOR = {Michiels, Wim and Niculescu, Silviu-Iulian},
  LOCATION = {Philadelphia},
  PUBLISHER = {Society for Industrial and Applied Mathematics},
  DATE = {2014},
  DOI = {10.1137/1.9781611973631},
  NUMBER = {27},
  SERIES = {Advances in Design and Control},
  SHORTSERIES = {Adv. Des. Control},
  SUBTITLE = {An Eigenvalue-Based Approach, second edition},
  TITLE = {Stability, Control, and Computation for Time-Delay Systems},
}

@BOOK{RoydenFitzpatrick2010,
  AUTHOR = {Royden, Halsey Lawrence and Fitzpatrick, Patrick M.},
  PUBLISHER = {Prentice Hall},
  DATE = {2010},
  EDITION = {4},
  TITLE = {Real Analysis},
}

@ARTICLE{SieberEngelborghsLuzyaninaSamaeyRoose2014,
  AUTHOR = {Sieber, Jan and Engelborghs, Koen and Luzyanina, Tatyana and Samaey, Giuseppe and Roose, Dirk},
  DATE = {2014},
  EPRINT = {1406.7144},
  EPRINTTYPE = {arXiv},
  JOURNALTITLE = {ArXiv e-prints},
  SUBTITLE = {Bifurcation analysis of delay differential equations},
  TITLE = {{DDE-BIFTOOL} manual},
}

@BOOK{Smith2011,
  AUTHOR = {Smith, Hal},
  LOCATION = {New York},
  PUBLISHER = {Springer},
  DATE = {2011},
  DOI = {10.1007/978-1-4419-7646-8},
  SERIES = {Texts in Applied Mathematics},
  SHORTSERIES = {Texts Appl. Math.},
  TITLE = {An Introduction to Delay Differential Equations with Applications to the Life Sciences},
}

@BOOK{Stepan1989,
  AUTHOR = {Stépán, Gábor},
  LOCATION = {Harlow},
  PUBLISHER = {Longman Scientific \& Technical},
  DATE = {1989},
  NUMBER = {210},
  SERIES = {Pitman Research Notes in Mathematics Series},
  SHORTSERIES = {Pitman Res. Notes Math. Ser.},
  TITLE = {Retarded Dynamical Systems: Stability and Characteristic Functions},
}

@BOOK{Trefethen2000,
  AUTHOR = {Trefethen, Lloyd Nicholas},
  LOCATION = {Philadelphia},
  PUBLISHER = {Society for Industrial and Applied Mathematics},
  DATE = {2000},
  DOI = {10.1137/1.9780898719598},
  SERIES = {Software, Environments and Tools},
  SHORTSERIES = {Software Environ. Tools},
  TITLE = {Spectral Methods in MATLAB},
}

@ARTICLE{Trefethen2008,
  AUTHOR = {Trefethen, Lloyd Nicholas},
  DATE = {2008},
  DOI = {10.1137/060659831},
  JOURNALTITLE = {SIAM Review},
  NUMBER = {1},
  PAGES = {67--87},
  SHORTJOURNAL = {SIAM Rev.},
  TITLE = {Is {G}auss Quadrature Better than {C}lenshaw--{C}urtis?},
  VOLUME = {50},
}
\par}

\end{document}